\documentclass{amsart}
\usepackage{amsfonts}
\usepackage[dvips]{graphicx}
\usepackage{epsfig} 
\newcommand{\Ha}{{\mathcal H}}
\newcommand{\Ka}{{\mathcal K}}
\newcommand{\jed}{{\mathbf 1}}
\newcommand{\C}{{\mathbb C}}
\newcommand{\R}{{\mathbb R}}
\newcommand{\E}{{\mathbb E}}
\newcommand{\N}{{\mathbb N}}
\newcommand{\M}{{\mathfrak M}}
\newcommand{\End}{{\rm End}}
\newcommand{\T}{{\mathbf T}}
\newcommand{\tr}{{\rm tr}}
\newcommand{\Tr}{{\rm Tr}}
\newcommand{\gwia}{^{\star}}
\newcommand{\ii}{{\mathbf i}}
\newcommand{\jj}{{\mathbf j}}
\newcommand{\kk}{{\mathbf k}}
\newcommand{\lll}{{\mathbf l}}
\newcommand{\equal}[2]{[#1=#2]}         
\newcommand{\equalA}[2]{\delta_{#1 #2}} 
\newcommand{\equalM}[2]{\delta_{#1 #2}} 
\newtheorem{definition}{Definition}
\newtheorem{theorem}{Theorem}
\newtheorem{remark}{Remark}
\newtheorem{proposition}{Proposition}
\newtheorem{lemma}{Lemma}

\begin{document}
\sloppy  
\title{Gaussian Random Matrix Models for $q$--deformed Gaussian Variables} 
\author{Piotr \'Sniady}
\address{Instytut Matematyczny, Uniwersytet Wroc\l{}awski, 
pl.~Grunwaldzki 2/4,\\ 
50-384 Wroclaw, Poland} \email{psnia@math.uni.wroc.pl}

\begin{abstract}
We construct a family of 
random matrix models for the $q$--deformed Gaussian random 
variables $G_\mu=a_\mu+a\gwia_\mu$ 
where the annihilation operators $a_\mu$ and creation operators $a\gwia_\nu$
fulfil the $q$--deformed commutation relation $a_\mu a\gwia_\nu-q 
a\gwia_\nu a_\mu=\Gamma_{\mu\nu}$, $\Gamma_{\mu\nu}$ is the covariance and
$0<q<1$ is a given number.
Important feature of considered random matrices is that
the joint distribution of their entries is Gaussian.
\end{abstract}

\maketitle 

\section{Introduction}
\subsection{The deformed Gaussian variables}
The $q$-deformed Gaussian random variables $G_\mu=a_\mu+a\gwia_\mu$, where
operators $a_\mu$ and their adjoints $a_\mu\gwia$  fulfil
deformed commutation relations \begin{equation} a_\mu a\gwia_\nu-q
a\gwia_\nu a_\mu=\Gamma_{\nu\mu} \jed \label{eq:qdeform}
\end{equation}
were introduced by Bourret and Frisch \cite{FB}. These operators act 
on a Hilbert space $\Ka$ which has a unital vector $\Omega$, called a vacuum, 
with the property that 
\begin{equation} \label{eq:anihilacja} a_\mu \Omega=0 \end{equation} 
for every value of the index $\mu$.

With the help of the vector $\Omega$ one can introduce a state $\tau$ on the 
algebra of operators acting on $\Ka$ as follows
$\tau(X)=\langle\Omega,X\Omega\rangle$. The state $\tau$ plays a role of
the non commutative expectation value.
From (\ref{eq:qdeform}) and (\ref{eq:anihilacja}) 
follows \cite{BS0} that for any $m\in\N$ and any indexes 
$\mu_1,\dots,\mu_{2m}$ we have that
\begin{equation} \label{eq:momenty1} \tau(G_{\mu_1} \cdots G_{\mu_{2m-1}})=0,
\end{equation}
\begin{equation} \label{eq:momenty2} \tau(G_{\mu_1} \cdots G_{\mu_{2m}})=
\sum_{\pi} q^{i(\pi)} \Gamma_{c_1 d_1} \cdots \Gamma_{c_m d_m},
\end{equation}
where the sum is taken over all pair partitions $\pi=\big\{\{c_1,d_1\},\dots,
\{c_m,d_m\}\big\}$ of the set $\{1,\dots, 2m\}$
and $i(\pi)$ is the number of crossings of the partition
$\pi$. For the reader's convenience we shall recall definitions of a pair 
partition and of its number of crossings in Sect.\ 
\ref{sec:definicjapartycji}.

From the quantum probability point of view all the information about 
non commutative random variables $G_\mu$ is encoded in their moments
$\tau(G_{\mu_1} \cdots G_{\mu_m})$ and therefore Eq.\ (\ref{eq:momenty1})
and (\ref{eq:momenty2}) can be treated as an alternative definition of 
$q$-deformed Gaussian variables $G_\mu$.

\subsubsection{Applications of deformed Gaussian variables}
Eq.\ (\ref{eq:momenty1}) and (\ref{eq:momenty2}) show that for $q=1$ 
operators $G_\mu$ have the 
same moments as classical Gaussian variables with the mean zero and the 
covariance $\Gamma_{\mu\nu}$, what should explain why do
we call $G_\mu$ deformed Gaussian variables. 
Eq.\ (\ref{eq:qdeform}) is for $q=1$ called the canonical 
(or bosonic) commutation relation.

For other special choices of the deformation parameter $q$ variables 
$G_\mu$ also have natural probabilistic interpretations \cite{FB}, namely
as increments of a dichotomic Markov process (for $q=-1$) or as Wigner's large 
random matrices (for $q=0$). 
Voiculescu \cite{V2}
has made a remarkable observation that for 
$q=0$ random variables $G_\mu$ are free semicircular elements (an analogue
of independent Gaussian variables in the free probability theory of Voiculescu 
\cite{V4,VDN}).  
Eq.\ (\ref{eq:qdeform}) is for $q=-1$ called
the canonical anticommutation relation (or fermionic relation) 
and for $q=0$ is called the free relation.

Therefore it was natural to expect that the relations (\ref{eq:momenty1}) and 
(\ref{eq:momenty2}) which 
are a simple generalisation of the three mentioned above: bosonic, fermionic 
and free cases
would give rise to interesting probabilistic objects. 

Indeed, it was observed by Bo\.zejko and Speicher \cite{BS0} that
a related to Eq.\ (\ref{eq:qdeform}) Brownian motion 
is a one component of an n-dimensional
Brownian motion which is invariant under the quantum 
group $S\sb{\sqrt{q}}U(n)$ of Woronowicz for $0<q<1$.

Another application of $q$-deformed Gaussian variables, this time 
as generalised quantum statistics was
proposed by Greenberg \cite{Greenberg} and Speicher \cite{Speicher}. 

The existence of operators $a_\mu$ and $a_\mu\gwia$ fulfilling deformed 
commutation relations 
(\ref{eq:qdeform}) was proven by Bo\.zejko and Speicher \cite{BS}. Later it was 
proven by Bo\.zejko, 
K\"ummerer, and Speicher \cite{BKS} 
that the von Neumann algebra generated by $q$--deformed Gaussian variables 
$G_1,G_2,\dots$ ($-1<q<1$) is a $II_1$ factor. There are today many open 
questions concerning these factors, particularly if they 
are different from the free group factors.

In this paper we present a natural probabilistic representation of
the $q$--deformed Gaussian variables for all $q\in [0,1]$ as 
some random matrices, what was
one of the open questions posed in  the paper \cite{FB}. A remarkable 
property of our model is that the joint distribution of entries of our 
matrices is Gaussian.

Recently a related problem of finding a random matrix model for so called 
$q$--deformed circular system was solved by Mingo and Nica \cite{MN}.

\subsubsection{The covariance $\Gamma_{\mu\nu}$}
Indexes $\mu,\nu$ are elements of a certain set $\M$. 
A necessary and sufficient condition
for operators $G_\mu$ to exist is that the function $\Gamma_{\mu\nu}$ is
positive definite \cite{BS}, i.e.\ 
$$\sum_{1\leq i,j\leq n} \alpha_i \alpha_j \Gamma_{\mu_i \mu_j}\geq 0$$
for all $\alpha_1,\dots,\alpha_n\in\R$ and $\mu_1,\dots,\mu_n\in\M$.
Typical examples of sets $\M$ and covariance functions are: 
\begin{itemize}
\item $\M=\N$ and $\Gamma_{i,j}=\delta_{ij}$. For $q=1$ we have that
$G_1,G_2,\dots$ is a sequence of independent, standard Gaussian variables
while for $q=0$ we have that $G_1,G_2,\dots$ is a sequence of free semicircular 
elements \cite{VDN}.
\item $\M=\R_+$ and $\Gamma_{t,s}=\mbox{min}(t,s)$. For $q=1$ we have that
$G_t$ is a Brownian motion, for $q=0$ we have that $G_t$ is a noncommutative 
stochastic process with 
free increments
\item $\M$ is a real Hilbert space and the covariance is defined by the scalar 
product: $\Gamma_{\phi\psi}=\langle \phi,\psi\rangle$. The case $\M={\mathcal 
L}^2(\R_+)$ is often used in the white noise calculus.
\end{itemize}

\subsubsection{The distribution of a deformed Gaussian variable}
A distribution of a random variable corresponding to the bounded
selfadjoint operator $G$ is a measure $\nu$ supported 
on the real line $\R$ such that
$$\tau(G^n)=\int x^n d\nu(x)$$
for all $n\in\N$.

It can be shown \cite{Sz} that the distribution $\nu_q$ of a $q$-deformed
Gaussian variable with the variance equal to $1$ is given by a measure
$\nu_q$ supported on the 
interval $\left[-\frac{2}{\sqrt{1-q}},\frac{2}{\sqrt{1-q}}\right]$
with a density
$$\nu_q(dx)=\frac{1}{\pi} \sqrt{1-q}\ \sin \theta \prod_{n=1}^{\infty}
(1-q^n) |1-q^n e^{2 i \theta}|^2 dx, $$
where $$x=\frac{2}{\sqrt{1-q}} \cos \theta $$ with $\theta\in [0,\pi]$.
  
\subsubsection{Canonical commutation relations and It\^o's formula} 
It is not merely an accident that for $q=1$ there is a correspondence between 
the commutation relation (\ref{eq:qdeform}) and Gaussian random variables.

If we consider a probability space generated by a Brownian motion $B(t)$ 
then every real random variable $X$ 
with a finite second moment can be uniquely expressed as a series of 
iterated It\^o integrals $$X=X^{(0)}+\sum_{i=1}^{\infty} 
\int_{0\leq t_1\leq \cdots \leq t_{i}<\infty} 
X^{(i)}(t_1,\dots,t_i) dB(t_1)  \cdots dB(t_i),$$ 
where $X^{(0)}\in\R$ and for each $i\in\N$ 
we have that $X^{(i)}:(\R_{+})^i\rightarrow \R$ is a symmetric function 
of its $i$ arguments. By the 
bosonic Fock space we call the space of sequences $(X^{(i)})_{i\geq 0}$ such 
that $X^{(i)}:(\R_{+})^i\rightarrow\R$ is a symmetric function. 

The Fock space carries a structure of a Hilbert space
with a scalar product 
$$\langle X,Y \rangle=\E[XY]=X^{(0)} Y^{(0)}+$$ 
$$\sum_{n\geq 1} \frac{1}{n!} \int_0^{\infty}\cdots \int_0^{\infty} 
X^{(n)}(t_1,\dots,t_n) Y^{(n)}(t_1,\dots,t_n) dt_1\cdots dt_n.$$

The Fock space representative of the random variable constantly equal to $1$ will 
be denoted by $\Omega$. We have $\Omega^{(0)}=1$,
$\Omega^{(n)}(t_1,\dots,t_n)=0$ for every $n\geq 1$. Note that for every random
variable we have
$$\E[X]=\E[X 1]=\langle \Omega,X \rangle.$$
 
It is often convenient to work  with the Fock space representations than with
random variables themselves. An interesting question is to determine the Fock 
space representation of a
product $[\int_0^\infty \phi(t) dB(t)] X$ if the Fock space representation of
$X$ is given and the integral is taken in the It\^o sense. 
It is a simple corollary from the It\^o's theorem that the
multiplication by $\int_0^\infty \phi(t) dB(t)$ can be expressed as a sum of
two operators acting on the Fock space: so called annihilation operator
$a_\phi$ and its adjoint, creation operator $a\gwia_\phi$: $$\int_0^\infty
\phi(t) dB(t)=a_\phi + a\gwia_\phi.$$
They have the following properties:  
\begin{equation} \label{eq:ccr} 
a_\phi a\gwia_\psi-a\gwia_\psi a_\phi= \langle \phi,\psi\rangle, \end{equation}
\begin{equation} a_\phi \Omega=0. 
\end{equation}

We see that the commutation relations (\ref{eq:ccr}) are equivalent to the 
statement of the It\^o's theorem. Conversely, postulating commutation relations 
(\ref{eq:qdeform}) is equivalent to saying that the non commutative stochastic 
process $G_t$ fulfils some non commutative It\^o's formula. Such a stochastic 
calculus for $q$--deformed operators was considered by the author \cite{Sn1}.

\subsection{Random matrices}
For a general introduction to random matrices and their applications in 
mathematics and physics
we refer to the monographs of Mehta \cite{M}, Hiai and Petz \cite{HP} 
and to overview articles 
\cite{Br,GMGW}.

There are essentially two kinds of questions concerning eigenvalues of a random matrix 
one can ask. Global questions are of the type: what is the asymptotic 
distribution of eigenvalues if the size of a matrix is large enough, 
while local questions concern, for example, the distribution
of the spacings between consecutive eigenvalues. The global questions are 
much easier to answer and free probability theory has provided powerful 
tools to answer
such questions for many random matrix model \cite{V3,V4,VDN,Sh}.

Recently random matrices were used as a powerful tool in the theory of 
$C\gwia$--algebras by Haagerup 
and Thorbjoersen \cite{HT}.

\subsection{Overview of the paper}
In Sect.\ \ref{sec:heureza} (which is independent of the rest of this article)
we present heuristic motivations 
for some matrix models considered in this article.

In Sect.\  \ref{sec:model} we introduce the notations and construct
a family we construct an auxiliary family of Gaussian random
matrices $R^{(N),A,\mu}$. 

In Sect.\ \ref{sec:dowod} we construct a central object of this paper,
namely a family of Gaussian random matrices $S^{(N),\mu}$. Since the
structure of matrices $S$ is a bit complicated, it is convenient to
think about them as some weighted sums of the auxiliary matrices
$R$, which have a much easier structure.

As the index $N$ tends to infinity, the size of our matrices grows
exponentially. We show that matrices $S$
asymptotically have the same expectation values as $q$--deformed Gaussian
random variables.

Our construction bases on the observation that for a finite dimensional 
Hilbert space $\Ha$
there are $2^N$ decompositions of a tensor power into two factors
$\Ha^{\otimes N}=\Ha^{\otimes k}\otimes \Ha^{\otimes (N-k)}$ which correspond
to ways of decomposing a set $\{1,\dots,N\}=A\cup (\{1,\dots,N\}\setminus A)$
into two subsets. The appropriate isomorphisms $j_A:\Ha^{\otimes N}\rightarrow
\Ha^{\otimes
|A|}\otimes\Ha^{\otimes (N-|A|)}$ give rise to isomorphisms of 
matrix algebras
$\tilde{j}_A:\End(\Ha^{\otimes |A|})\otimes \End(\Ha^{\otimes (N-|A|)})
\rightarrow \End(\Ha^{\otimes N})$.

Each auxiliary matrix $R^{(N),A,\mu}$ is obtained by embedding
a small standard hermitian random matrix 
$R_1^{(N),\mu,A}\in\End (\Ha^{|A|} ) $ into a bigger algebra $\End (\Ha^N)$.
The embedding is implemented by the isomorphism $\tilde{j}_A$.
We have that
$$S^{(N),\mu}=\sum_{A\subseteq\{1,\dots,N\} } \sigma^{(N)}_A R^{(N),A,\mu}=
\sum_{A\subseteq\{1,\dots,N\} } \sigma^{(N)}_A
\tilde{j}_A [R_1^{(N),A,\mu}\otimes \jed], $$
where $\sigma^{(N)}_A$ is a certain weight function.
 
It turns out that the the
commutation properties of two random matrices 
$R^{(N),A_i,\mu}$ corresponding to two decompositions given by
sets $A_1, A_2$ depend on the number of common elements of $A_1$ and $A_2$.
For example, if $A_1\cap A_2=\emptyset$ then these two matrices commute and,
informally speaking, the more elements $A_1$ and $A_2$ have in common, the
more they behave like a pair of free random variables. Therefore the
expectation value of a product of many  matrices $R$ can be evaluated from
the number of elements of $A_i\cap A_j$ and in this way we are able
to calculate the mixed moments of matrices $S$.

The choice of a normed weight $\sigma^{(N)}_A$ is equivalent to a choice of a
probabilistic measure $\rho$ on the set of all subsets of $\{1,\dots,N\}$
defined such that the measure of a singleton $\{A\}$ is equal to 
$(\sigma^{(N)}_A)^2$ for
any $A\subseteq\{1,\dots,N\}$.
If this measure is, loosely speaking, concentrated on the sets of order
$c\sqrt{N}$ then if $A$ and $B$ are independent random variables with
distribution given by the measure $\rho$ then $|A\cap B|$ is asymptotically
Poisson distributed with the parameter $\lambda=c^2$.  For appropriate
choice of $c$ we are able in this way to obtain a random matrix model for
$q$--deformed Gaussian random variables for all $0\leq q\leq 1$.

In Sect.\ \ref{sec:almostsurely} we show that our matrices converge 
to $q$-deformed Gaussian variables not only in the sense of expectation
values mixed moments, but that that mix moments converge almost surely. 

Sect.\ \ref{sec:technicznedowody} is devoted to proofs of technical 
lemmas.

%
The random matrices $S$ considered in this article are complex hermitian. However there
are no difficulties to extend these results to real symmetric or symplectic hermitian.

This paper shows that the $q$--deformed probability theory, which was
regarded until today as purely abstract and algebraic, in fact has natural
probabilistic models just like the free probability theory has.

\section{Heuristics of 
Random Commutation Relations and Random Gaussian Matrices}
\label{sec:heureza}
This section an an independent part of the article and notations used 
here will not be used in the subsequent considerations. 
We have to warn the reader that this section is very informal, however by such
informal considerations it is much easier to get an insight to the nature of
the problem.

\subsection{$q$--deformed Gaussian variables}
Our motivations how to find random matrices which asymptotically behave 
like $q$--deformed Gaussian variables were inspired by a careful study of 
the article of Speicher \cite{Sp}. 

In this paper he shows a certain 
non commutative central limit theorem that if a suitably normalised family 
of centered non commutative random variables $K_1, K_2, \dots$ 
has a property that 
each pair of them either commutes (with 
probability $p$) or anticommutes (with probability $1-p$) and if for each pair
the choice of 
one these possibilities is made independently, then the distribution
of a normalized mean $\frac{K_1+\dots+K_n}{\sqrt{n}}$ converges 
(as $n\rightarrow\infty$)
to the distribution of a $q$--deformed Gaussian 
random variable with $q=2p-1$. 

Now we shall construct a family of random (non Gaussian) 
matrices which almost fulfils the assumptions of the Speicher's theorem.

Let us fix a real number $0<q<1$. For any $N\in\N$ let us consider a 
family of $2^N\times 2^N$ matrices  
$$K_s=K_s^1\otimes \cdots \otimes K_s^N$$
indexed by $s\in\N$, 
where for each $s\in\N$ and $1\leq n\leq N$ we have that $K_s^n$ is a $2\times 
2$ matrix chosen randomly according to the following table:

\begin{center}
\begin{tabular}{c|c|c|c|c|c|c|c|c}
matrix & $\sigma_0$ & $-\sigma_0$ & $\sigma_1$ & $-\sigma_1$ 
& $\sigma_2$ & $-\sigma_2$ & $\sigma_3$ & $-\sigma_3$ \\
\hline
probability & $\frac{1-3r}{2}$ & $\frac{1-3r}{2}$ & $\frac{r}{2}$ & 
 $\frac{r}{2}$ &  $\frac{r}{2}$ &  $\frac{r}{2}$ &  $\frac{r}{2}$ &  
$\frac{r}{2}$ 
\end{tabular},
\end{center}

where $0\leq r\leq \frac{1}{3}$ is a real number to be specified later and
$$\sigma_0=\left[ \begin{array}{cc} 1 & 0 \\ 0 & 1 \end{array} \right],\quad 
\sigma_1=\left[ \begin{array}{cc} 0 & 1 \\ 1 & 0 \end{array} \right],  \quad
\sigma_2=\left[ \begin{array}{cc} 0 & -i \\ i & 0 \end{array} \right], \quad
\sigma_3=\left[ \begin{array}{cc} 1 & 0 \\ 0 & -1 \end{array} \right] $$ 
are Pauli matrices. The random choices of matrices $K_s^n$ 
should be made independently.

These eight hermitian matrices $\pm \sigma_i$ have a property that each pair 
of them either commutes or anticommutes. It is a simple calculation that
the probability that two independent $2\times 2$ matrices (chosen according to
the above table) anticommute is equal to $6r^2$.

Therefore all matrices $K_s$ are hermitian and furthermore each pair 
$K_s$ and $K_t$ ($s\neq t$) either commutes (if the number of indexes $n$
such that $K_s^n$ anticommutes with $K_t^n$ is even) or anticommutes
(if this number is odd). We see that 
the probability of the event that $K_s$ commutes with $K_t$ is equal to
$$p=\sum_{0\leq m \leq \frac{N}{2}} {N \choose 2m} (1-6r^2)^{N-2m} 
(6r^2)^{2m}, $$
while the probability that $K_s$ anticommutes with $K_t$ is equal to
$$1-p=\sum_{0\leq m \leq \frac{N-1}{2}} {N \choose 2m+1} (1-6r^2)^{N-(2m+1)} 
(6r^2)^{2m+1}.$$
The difference of these two probabilities is equal to $2p-1=(1-12r^2)^N$,

We would like to apply now the Speicher's theorem to the family $(K_s)$ 
by choosing
$r$ as a function of $N$ such that $q=[1-12r_N^2]^N$. For large 
$N$ one can take the approximate value $r_N=\sqrt{-\frac{\ln q}{12 N}}$.

Unfortunately the Speicher's theorem cannot be applied directly
since it turns out that the events ``$K_s$ commutes with $K_t$'' 
are not independent for different pairs $\{s,t\}$. 
However, if $N$ tends to infinity it can be justified that they are, 
loosely speaking, more and more independent.

By classical central limit theorem we have that the limit distribution of 
$K=\frac{1}{\sqrt{n}}(K_1+\cdots+K_n)$ (regarded as a classical random 
variable with
values in a vector space of matrices) is Gaussian. In order to characterise
this distribution uniquely we have  to give the mean and the covariance of 
entries. The mean value of each entry 
of the matrix $K$ is equal to $0$ and the covariance factorises as follows
$$\E[K_{(i_1,\dots,i_N),(j_1,\dots,j_N)} K_{(k_1,\dots,k_N),(l_1,\dots,l_N)}]=$$
$$=\E[K_{(i_1,\dots,i_N),(j_1,\dots,j_N)} 
\overline{K_{(l_1,\dots,l_N),(k_1,\dots,k_N)}}]=$$
$$\E[K^1_{i_1,j_1} K^1_{k_1,l_1}] \cdots \E[K^N_{i_N,j_N} K^N_{k_N,l_N}].$$
Above we have used the fact that due to the factorisation 
$M_{2^N}=M_2\otimes\cdots\otimes M_2$ we can parametrise the rows and columns
of a $2^N\times 2^N$ matrix by sequences $(i_1,\dots, i_N)$, where
$0\leq i_1,\dots, i_N\leq 1$. 

A simple calculations shows that the covariance of $2\times 2$ matrices 
given by the table above is equal to
$$\E[K^n_{i_n,j_n} K^n_{k_n,l_n}]=(1-4r) \delta_{i_n,j_n} \delta_{k_n,l_n}+
4r \frac{1}{2} \delta_{i_n,l_n} \delta_{j_n,k_n}.$$

It turns out that this is exactly the covariance we will 
considered in Proposition \ref{prop:zbieznosc} for the special case of $d=2$, 
see (\ref{eq:przykladkowariancji}).


\section{Random Matrix Model for Deformed Gaussian Variables}
\label{sec:model}
\label{sec:definicjapartycji}
\subsection{Pair partitions}
\begin{definition} A pair partition of a given finite set $M$ is any
decomposition of $M$ into a family $\pi=\bigl\{ \{c_1,d_1\},\dots,\{c_m,d_m\}
\bigr\}$ of disjoint sets each having exactly two elements:  $$c_i\neq
d_i,\qquad \mbox{for }1\leq i\leq m, $$
$$\{c_i,d_i\}\cap\{c_j,d_j\}=\emptyset, \qquad \mbox{for }i\neq j,$$
$$\{c_1,d_1\}\cup\cdots\cup\{c_m,d_m\}=M.$$

The sets $\{c_1,d_1\},\dots,\{c_m,d_m\}$ are called lines of the pair
partition $\pi$.

If $M$ is an ordered set we say that two distinct lines $\{a,b\},\{u,v\}$, 
$a<b$, $u<v$ cross if
$a<u<b<v$ or $u<a<v<b$.

For a given pair partition $\pi$ we will denote by $i(\pi)$ the number of
crossings in $\pi$, i.e. number of all unordered pairs of lines 
$\{a,b\},\{u,v\}\in\pi$ such that the lines $\{a,b\}, \{u,v\}$ cross. 
\end{definition} 
{\bf Example.} There is only one pair partition of a two 
element set $\{1,2\}$ and there are three pair partitions of a four element 
set $\{1,2,3,4\}$, namely $\bigl\{ \{1,2\},\{3,4\} \bigr\}$,  
$\bigl\{ \{1,4\},\{2,3\} \bigr\}$, $\bigl\{\{1,3\},\{2,4\} \bigr\}$.


We have $i\left(\bigl\{ \{1,2\},\{3,4\} \bigr\}\right)=0$, 
$i\left(\bigl\{\{1,3\},\{2,4\} \bigr\}\right)=1$.

Sets having an odd number of elements do not have any pair partitions at all.

\subsection{Notations} Let us fix a natural number $d\geq 2$. For any 
$N\in\N$ we define $\Ha_r=\C^d$ for $1\leq r\leq N$ and 
$\Ha^{(N)}=\bigotimes_{1\leq r\leq N} \Ha_r$. In the following we shall 
often omit the index $(N)$ standing at various objects, however we have to 
remember about their dependence on $N$.

If $f_0,\dots,f_{d-1}$ is an orthonormal basis of $\C^d$ then
$e_{\ii}=f_{i_1}\otimes\cdots\otimes f_{i_N}$ is an orthonormal basis of
$\Ha^{(N)}$, where $\ii=(i_1,\dots,i_N)$ and $0\leq i_1,\dots,i_N\leq d-1$.

Here and in the following by bold letters $\ii,\jj,\kk,\dots$ we shall denote
variables which index the basis of $\Ha^{(N)}$ and always we have 
$\ii=(i_1,\dots,i_N)$, $\kk=(k_1,\dots,k_N)$, etc. 

For any set $A$ where  $A=\{a_{1},\dots,a_k\}\subseteq\{1,\dots,N\}$,
$a_{1}<\cdots<a_k$ and $A'=\{1,2,\dots,N\}\setminus A=\{b_1,\dots,b_{N-k}
\}$, $b_1<\cdots<b_{N-k}$ we consider Hilbert spaces
$\Ha^A=\bigotimes_{r\in A} \Ha_r$ and $\Ha^{A'}=\bigotimes_{r\in A'} \Ha_r$
and an isomorphism $j_A:\Ha^{(N)}\rightarrow \Ha^A\otimes\Ha^{A'}$ given by
a grouping of factors:
$$v_1\otimes\cdots\otimes v_N\mapsto (v_{a_{1}}\otimes\cdots\otimes
v_{a_k})\otimes (v_{b_1}\otimes\cdots\otimes v_{b_{N-k}}).$$ This
isomorphism induces an isomorphism of matrices
$\tilde{j}_A:\End(\Ha^A)\otimes\End(\Ha^{A'})\rightarrow\End(\Ha^{(N)})$.

In the following we shall denote by $\tr$ the normalised trace on
$\End(\Ha^{(N)})$ defined by $\tr(M)=\frac{1}{d^N} \Tr(M)$, where $\Tr$
denotes the standard trace.


\subsection{Iverson's notation} 
In the following we shall use sometimes Iverson's notation \cite{GKP} as an 
alternative to the Kronecker's notation:
$$[x=y]=\delta_{xy}=\left\{ \begin{array}{cl} 0 & \quad \mbox{if } x\neq y \\
1 & \quad \mbox{if } x=y \end{array} \right. .$$
Of course the Iverson's symbol $[x=y]$ is always equal to the Kronecker's
delta $\delta_{xy}$ but it has some typographic advantages if $x$ and $y$
are complicated expressions with many upper and lower indexes.

\subsection{Random matrices $R$}
\label{sec:losowe}
\begin{definition}  \label{def:matrix}
If $V$ is a finite dimensional Hilbert space with an orthonormal basis
$e_1,\dots, e_{{\rm dim\ } V}$  then a hermitian standard
random matrix over $V$ is a random variable $M$ with values in
$\End(V)$ such that the
joint distribution of the complex matrix coefficients 
$M_{ij}=\langle e_i, M e_j \rangle$
is Gaussian, $M_{ij}=\overline{M_{ji}}$, all $M_{ij}$ have mean zero and
the covariance is given by
\begin{equation} \E [M_{ij} M_{kl}]=\E [M_{ij} \overline{M_{lk}}]=
\frac{1}{{\rm dim\ }V} \equal{i}{l} \equal{j}{k}.
\label{eq:kowariancja0} \end{equation}

Alternatively one can define a hermitian standard random matrix
$(M_{ij})$ by saying that the following random variables: $M_{ii}$ for all 
indexes $i$, $\Re M_{ij}$, $\Im
M_{ij}$ for $i<j$ are independent real Gaussian variables with
$\E(M_{ij})=0$ for all indexes $i,j$ and $\E M_{ii}^2=\frac{1}{{\rm dim\ }V}$
for all values of index $i$ and $\E (\Re M_{ij})^2=\E (\Im
M_{ij})^2=\frac{1}{2 {\rm \ dim\ }V}$ for all $i<j$. The entries
$M_{ij}$, where $i>j$ are defined by the hermitianity condition
$M_{ij}=\overline{M_{ji}}$.

One can show that both definitions do not depend on the choice of the 
orthonormal basis of $V$.
\end{definition} 

For each $A\subseteq\{1,\dots,N\}$ let us consider a family of hermitian
standard random matrices $R^{(N),A,\mu}_1\in \End (\Ha_A)$ indexed by
$\mu\in\M$ such that the entries of different
matrices are independent. We define a family of random matrices 
$R^{(N),A,\mu}$ by
$$R^{(N),A,\mu}=\tilde{j}_A (R^{(N),A,\mu}_1 
\otimes\jed_{\Ha^{A'}})\in\End(\Ha^{(N)})$$
where $\jed_{\Ha^{A'}}:\Ha^{A'}\rightarrow \Ha^{A'}$ denotes the identity
operator.

Intuitively speaking, a matrix $R^{(N),A,\mu}$ consists of $d^{N-|A|}$
copies of a $d^{|A|}\times d^{|A|}$ standard hermitian random matrix.

As one can see, matrices $R^{A,\mu}$ are hermitian and
the joint distribution of their entries is Gaussian, but different entries need
not to be independent. We have:
$$R^{A,\mu}_{\ii\jj}=\overline{R^{A,\mu}_{\jj\ii}},
\qquad \E [R^{A,\mu}_{\ii\jj}]=0, $$
and from (\ref{eq:kowariancja0}) it follows that
\begin{equation} \E [R^{A,\mu}_{\ii\jj} R^{B,\nu}_{\kk\lll}]=\E
[R^{A,\mu}_{\ii\jj} 
\overline{R^{B,\nu}_{\lll\kk}}]=\label{eq:kowariancja1} \end{equation} 
$$= \equalA{A}{B} \equalM{\mu}{\nu} \left(\prod_{r\in A} 
\frac{\equal{i_r}{l_r} \equal{j_r}{k_r}}{d} \right) \left(\prod_{r\in A'} 
\equal{i_r}{j_r} \equal{k_r}{l_r}\right). $$

\subsection{Tensors $T$}
The formula (\ref{eq:kowariancja1}) can be written shorter if we introduce for 
all $A\subseteq\{1,\dots,N\}$ and $1\leq r\leq N$ tensors $T_{ij,kl}^{A,r}$ as 
follows:
$$T_{ij,kl}^{A,r}=
\left\{\begin{array}{cl} \frac{1}{d} \equal{i}{l} \equal{j}{k} &\quad \mbox{if 
} r\in A \\ 
\equal{i}{j} \equal{k}{l} &\quad \mbox{if } r\not\in A \end{array} 
\right. .$$
We define  
\begin{equation} \T_{\ii\jj,\kk\lll}^{A}=\prod_r T_{i_r j_r,k_r l_r}^{A,r},
\label{eq:tensort1} \end{equation}
what with a small abuse of notation can be written as
\begin{equation} \label{eq:tensort2} \T^{A}=T^{A,1}\otimes \cdots\otimes 
T^{A,N}.
\end{equation}

Then (\ref{eq:kowariancja1}) can be written as
\begin{equation} \E [R^{A,\mu}_{\ii\jj} R^{B,\nu}_{\kk\lll}]=\E
[R^{A,\mu}_{\ii\jj} 
\overline{R^{B,\nu}_{\lll\kk}}]=\equalA{A}{B} \equalM{\mu}{\nu} 
\T^A_{\ii\jj,\kk\lll}.\label{eq:kowariancja1a} \end{equation} 
\subsection{Examples}
First of all note that for the trivial case $d=1$ all Hilbert spaces are 
one dimensional 
and all random matrices $R^{A,\mu}$ are in fact scalar random Gaussian 
variables. 

Furthermore, the random matrix $R^{(N),\mu,\emptyset}$ is simply a 
scalar real random Gaussian variable multiplied by an identity matrix. 
The random matrix $R^{(N),\mu,\{1,\dots,N\} }$ is a hermitian 
standard random matrix from Definition \ref{def:matrix}.

There is a correspondence between sequences $\ii=(i_1,\dots,i_N)$ such that
$0\leq i_1,\dots,i_N \leq d-1$ and the set of integer numbers
$\{0,1,\dots,d^N-1\}$ given by the digit representation of natural
numbers in the system with base $d$:
$$\ii=(i_1,\dots,i_N)\mapsto i_1+d i_2+\cdots+d^{N-1} i_{N}.$$ 
Therefore we can introduce an orthonormal basis $g_0,\dots,g_{d^N-1}$ of 
$\Ha^{(N)}$ indexed by integer numbers:
 $$g_{i_1+d i_2+\cdots+d^{N-1} i_N}=f_{i_1}\otimes\cdots\otimes 
f_{i_N}= e_{(i_1,\dots,i_N)},$$
where $0\leq i_1,\dots,i_N\leq d-1$. In the following, if we want to write 
an endomorphism $M\in\End(\Ha^{(N)})$ as a matrix $(M_{ij})_{0\leq i,j\leq 
d^N-1}$ we shall do it in the basis $(g_i)_{0\leq i\leq d^N-1}$.

For $d=2$ and $N=2$ the matrices $R^{(N),\mu,A}$ are of the
following form:
$$R^{\{1\} }= \left[ \begin{array}{cccc} {a_{00}} & {a_{01}} & 0 & 0 \cr
{a_{10}} & {a_{11}} & 0 & 0 \cr 0 & 0 &
{a_{00}} & {a_{01}} \cr 0 & 0 & {a_{10}} & {a_{11}} \end{array} \right], $$
$$R^{\{2 \} }= \left[ \begin{array}{cccc} {b_{00}} & 0 & {b_{01}} & 0 \cr 0 &
{b_{00}} & 0 & {b_{01}} \cr {b_{10}} & 0 &
{b_{11}} & 0 \cr 0 & {b_{10}} & 0 & {b_{11}} \end{array} \right], $$
where $\left[ \begin{array}{cc} a_{00} & a_{01} \\ a_{10} & a_{11} \end{array}
\right]$ and $\left[ \begin{array}{cc} b_{00} & b_{01} \\ b_{10} & b_{11}
\end{array} \right]$ are standard
hermitian random matrices. Entries of the first matrix are by definition
independent of the entries of the second matrix.
The index $\mu$ was omitted, however it should be
understood that for different values $\mu$ the entries of matrices are
independent.

For $d=2$ and $N=3$ we have:
$$R^{\{1\} }=\left[ \begin{array}{cccccccc} {c_{00}} & {c_{01}} & 0 & 0 & 0 & 0
& 0 & 0
\cr {c_{10}} & {c_{11}} &
0 & 0 & 0 & 0 & 0 & 0 \cr 0 & 0 & {c_{00}} & {c_{01}} & 0 & 0 & 0 & 0 \cr 0 & 0
& {c_{10}} & {c_{11}} &
0 & 0 & 0 & 0 \cr 0 & 0 & 0 & 0 & {c_{00}} & {c_{01}} & 0 & 0 \cr 0 & 0 & 0 & 0
& {c_{10}} & {c_{11}} &
0 & 0 \cr 0 & 0 & 0 & 0 & 0 & 0 & {c_{00}} & {c_{01}} \cr 0 & 0 & 0 & 0 & 0 & 0
& {c_{10}} & {c_{11}}
\end{array} \right], $$
$$R^{\{2\} }=\left[ \begin{array}{cccccccc} {d_{00}} & 0 & {d_{01}} & 0 & 0 & 0
& 0
& 0 \cr 0 & {d_{00}} & 0 &
{d_{01}} & 0 & 0 & 0 & 0 \cr {d_{10}} & 0 & {d_{11}} & 0 & 0 & 0 & 0 & 0
\cr 0 & {d_{10}} & 0 &
{d_{11}} & 0 & 0 & 0 & 0 \cr 0 & 0 & 0 & 0 & {d_{00}} & 0 & {d_{01}} & 0
\cr 0 & 0 & 0 & 0 & 0 &
{d_{00}} & 0 & {d_{01}} \cr 0 & 0 & 0 & 0 & {d_{10}} & 0 & {d_{11}} & 0 \cr 0 &
0 & 0 & 0 & 0 &
{d_{10}} & 0 & {d_{11}} \end{array} \right], $$
$$R^{\{3\} }=\left[\begin{array}{cccccccc} {e_{00}} & 0 & 0 & 0 & {e_{01}} & 0 &
0 & 0 \cr
0 &
{e_{00}} & 0 & 0 & 0 & {e_{01}} & 0 & 0 \cr 0 & 0 & {e_{00}} & 0 & 0 & 0 &
{e_{01}} & 0 \cr 0 & 0 & 0 &
{e_{00}} & 0 & 0 & 0 & {e_{01}} \cr {e_{10}} & 0 & 0 & 0 & {e_{11}} & 0 & 0 & 0
\cr 0 & {e_{10}} & 0 &
0 & 0 & {e_{11}} & 0 & 0 \cr 0 & 0 & {e_{10}} & 0 & 0 & 0 & {e_{11}} & 0 \cr 0
& 0 & 0 & {e_{10}} & 0 &
0 & 0 & {e_{11}}   \end{array} \right], $$
where again $(c_{pq})_{0\leq p,q\leq 1},(d_{pq})_{0\leq p,q\leq
1},(e_{pq})_{0\leq p,q\leq 1}$ are standard hermitian random matrices.

Furthermore
$$R^{\{1,2\} }=\left[ \begin{array}{cccccccc} {f_{00}} & {f_{01}} & 
{f_{02}} & {f_{03}} & 0 & 0 & 0 & 0 \cr {f_{10}} & {f_{11}} & {f_{12}} & 
{f_{13}} & 0 & 0 & 0 & 0 \cr {f_{20}} & {f_{21}} & {f_{22}} & {f_{23}} & 0 & 
0 & 0 & 0 \cr {f_{30}} & {f_{31}} & {f_{32}} & {f_{33}} & 0 & 0 & 0 & 0 \cr 
0 & 0 & 0 & 0 & {f_{00}} & {f_{01}} & {f_{02}} & {f_{03}} \cr 0 & 0 & 0 & 0 
& {f_{10}} & {f_{11}} & {f_{12}} & {f_{13}} \cr 0 & 0 & 0 & 0 & {f_{20}} & 
{f_{21}} & {f_{22}} & {f_{23}} \cr 0 & 0 & 0 & 0 & {f_{30}} & {f_{31}} & 
{f_{32}} & {f_{33}} \end{array} \right], $$ 
$$R^{\{1,3\}} =\left[ \begin{array}{cccccccc} {g_{00}} & {g_{01}} & 0 & 0 & 
{g_{02}} & {g_{03}} & 0 & 0 \cr {g_{10}} & {g_{11}} & 0 & 0 & {g_{12}} & 
{g_{13}} & 0 & 0 \cr 0 & 0 & {g_{00}} & {g_{01}} & 0 & 0 & {g_{02}} &
{g_{03}} \cr 0 & 0 & {g_{10}} & {g_{11}} & 0 & 0 & {g_{12}} & {g_{13}} \cr 
{g_{20}} & {g_{21}} & 0 & 0 & {g_{22}} & {g_{23}} & 0 & 0 \cr {g_{30}} & 
{g_{31}} & 0 & 0 & {g_{32}} & {g_{33}} & 0 & 0 \cr 0 & 0 & {g_{20}} & 
{g_{21}} & 0 & 0 & {g_{22}} & {g_{23}} \cr 0 & 0 & {g_{30}} & {g_{31}} & 0 & 
0 & {g_{32}} & {g_{33}} \end{array}\right], $$
$$R^{\{2,3\}}= \left[ \begin{array}{cccccccc} {h_{00}} & 0 & {h_{01}} & 0 &
{h_{02}} & 0 & {h_{03}} & 0 \cr 0 & {h_{00}} & 0 & {h_{01}} & 0 & {h_{02}} & 
0 & {h_{03}} \cr {h_{10}} & 0 & {h_{11}} & 0 & {h_{12}} & 0 & {h_{13}} & 0 
\cr 0 & {h_{10}} & 0 & {h_{11}} & 0 & {h_{12}} & 0 & {h_{13}} \cr {h_{20}} & 
0 & {h_{21}} & 0 & {h_{22}} & 0 & {h_{23}} & 0 \cr 0 & {h_{20}} & 0 & 
{h_{21}} & 0 & {h_{22}} & 0 & {h_{23}} \cr {h_{30}} & 0 & {h_{31}} & 0 & 
{h_{32}} & 0 & {h_{33}} & 0 \cr 0 & {h_{30}} & 0 & {h_{31}} & 0 &
{h_{32}} & 0 & {h_{33}} \end{array} \right], $$
where $(f_{pq})_{0\leq p,q\leq 3},(g_{pq})_{0\leq p,q\leq 3},(h_{pq})_{0\leq
p,q\leq 3}$ are standard hermitian random matrices.

\subsection{The case of a general covariance $\Gamma_{\mu \nu}$}
By a small change of definition of the matrices $R$ we obtain a more
general case.

Let $\Gamma_{\mu\nu}$ be a real positive definite function. 
For every $A\subseteq\{1,\dots,N\}$ we consider a family of random 
(non hermitian)
matrices $R^{(N),A,\mu}_0\in\End(\Ha_A)$ such that for each pair of indexes
$\ii,\jj$ we have that
the joint distribution of  $(R^{(N),A,\mu}_0)_{\ii\jj}$  $\mu\in\M$ is
Gaussian, $$\E(R^{(N),A,\mu}_0)_{\ii\jj}=0,$$ the covariance of real and
imaginary parts are defined by the function $\Gamma$:
  $$\E[
\Re(R^{(N),A,\mu}_0)_{\ii\jj} \Re(R^{(N),A,\nu}_0)_{\ii\jj}]= \E[
\Im(R^{(N),A,\mu}_0)_{\ii\jj} \Im(R^{(N),A,\nu}_0)_{\ii\jj}]=\frac{1}{2
d^{|A|}} \Gamma_{\mu\nu} $$
and the real and imaginary parts are independent:
$$\E[ \Re(R^{(N),A,\mu}_0)_{\ii\jj} \Im(R^{(N),A,\nu}_0)_{\ii\jj}]=0. $$

For different choices of sets $A$ or a pair of indexes $(\ii,\jj)$ the
random variables $(R^{(N),A,\mu}_0)_{\ii\jj}$ should be independent.

We define hermitian random matrices 
$$R^{(N),A,\mu}_1=R^{(N),A,\mu}_0+(R^{(N),A,\mu}_0)^{\star}.$$ 
Note that for the simplest choice
of a positive definite function $\Gamma_{\mu\nu}=\equalM{\mu}{\nu}$ this
definition of random matrices $R^{(N),A,\mu}_1$ coincides with the
definition from Subsect.\ \ref{sec:losowe}.

Similarly as in Subsect.\ \ref{sec:losowe} we define
$$R^{(N),A,\mu}=\tilde{j}_A (R^{(N),A,\mu}_1 \otimes\jed_{\Ha_{A'}}).$$
The joint distribution of entries of hermitian matrices $R^{(N),A,\mu}$ is
Gaussian and
$$\E[R^{A,\mu}]_{\ii\jj}=0,  $$
\begin{equation} \E [R^{A,\mu}_{\ii\jj} R^{B,\nu}_{\kk\lll}]=\E
[R^{A,\mu}_{\ii\jj} \overline{R^{B,\nu}_{\lll\kk}}]=\equalA{A}{B}\Gamma_{\mu 
\nu} \T^A_{\ii\jj,\kk\lll}. \label{eq:kowariancja}
\end{equation} 

\section{The Main Theorem}
\label{sec:dowod}
We define a family of random matrices indexed by $\mu\in\M$
$$ S^{(N),\mu}=\sum_{A\subseteq\{1,\dots,N\} } \sigma^{(N)}_A R^{(N),A,\mu} $$
where $\sigma^{(N)}$ is a real-valued function on the set of all subsets
of $\{1,\dots,N\}$.

Matrices $S^{(N)}$ are hermitian and the joint distribution of their 
entries is Gaussian. Alternatively one can define these matrices by giving
the mean and the covariance of the entries: we have 
\begin{equation} \E[S^{\mu}_{\ii\jj}]=0, \end{equation}
\begin{equation} \E[S^{\mu}_{\ii\jj} S^{\nu}_{\kk\lll}] =
\E[S^{\mu}_{\ii\jj} \overline{S^{\nu}_{\lll\kk}}]=
\Gamma_{\mu\nu} \sum_A (\sigma_A)^2 \T^A_{\ii\jj,\kk\lll}.
\label{eq:kowariancjas}
\end{equation}

In the following theorem we show conditions which the sequence of functions
$(\sigma^{(N)})$ need to fulfil. 
Since these conditions may seem quite disgusting, 
we would like to give some hope to the reader by pointing to the equation 
(\ref{eq:przykladkowariancji}), which gives a simple example of a 
covariance function fulfilling all assumptions.

\begin{theorem} \label{th:theo}
If for each $N\in\N$ we have that $\sigma^{(N)}$ is a real--valued 
function on the set of all subsets of $\{1,\dots,N\}$ such that:
\begin{enumerate}
\item \label{zal:z1} (normalisation) for each $N\in\N$ we have
$$\sum_{A\subseteq\{1,\dots,N\} } (\sigma^{(N)}_A)^2=1,$$
\item \label{zal:z2} (triple coincidations are rare)
$$\lim_{N\rightarrow\infty} \sum_{{A_1,A_2,A_3\subseteq\{1,\dots,N\}}\atop{\ 
A_1\cap A_2\cap A_3\neq\emptyset}} (\sigma^{(N)}_{A_1})^2 
(\sigma^{(N)}_{A_2})^2 (\sigma^{(N)}_{A_3})^2=0,$$
\item \label{zal:z3} (distribution of coincidations)
there exists a sequence $(p_i)_{i\geq 0}$ of nonnegative real 
numbers such that $\sum_{i\geq 0} p_i=1$ and for
any $k\in\N$ and any nonnegative integer numbers $n_{ij}$, $1\leq i<j\leq 
k$ we have
$$\lim_{N\rightarrow\infty} \sum_{{A_1,\dots,A_k\subseteq\{1,\dots,N\}}\atop{
|A_i\cap A_j|=n_{ij} {\rm \, for\: any\, } 1\leq i<j\leq k}} 
(\sigma^{(N)}_{A_1})^2 \cdots (\sigma^{(N)}_{A_k})^2=\prod_{1\leq i<j\leq k} 
p_{n_{ij}}, $$ 
\item \label{zal:z4} for each $n\in\N$  
$$\lim_{N\rightarrow\infty} \sum_{A_1,\dots,A_n\subseteq\{1,\dots,N\} } 
(\sigma^{(N)}_{A_1})^2 \cdots 
(\sigma^{(N)}_{A_n})^2 \frac{1}{d^{2|A_1\setminus (A_2\cup \cdots\cup 
A_n)|}}=0.$$
%
\end{enumerate}
Then for $q=\sum_{i=0}^{\infty} p_i \frac{1}{d^{2i}} $ we have 
$$\lim_{N\rightarrow\infty} \E[\tr\ S^{(N),\mu_1} \cdots S^{(N),\mu_n}] =
\sum_\pi q^{i(\pi)} \Gamma_{\mu_{c_1} \mu_{d_1}} \cdots \Gamma_{\mu_{c_m}
\mu_{d_m}}= 
\tau[G_{\mu_1} \cdots G_{\mu_n}]$$
for every $n\in\N$ and $\mu_1,\dots,\mu_n\in\M$,
where $G_{\mu_1},\dots,G_{\mu_n}$ are $q$-deformed Gaussian variables with covariance
$\Gamma$.
\end{theorem} 
Before we prove this theorem we would like to make some remarks and
state auxiliary lemmas.

\begin{remark} For a given function $\sigma^{(N)}$ 
we define a measure $\rho^{(N)}$ on the set of all subsets of 
$\{1,\dots,N\}$ by assigning to set $A$ the weight $(\sigma^{(N)}_A)^2$. 
Then the first three assumptions of the theorem can be reformulated as 
follows: \begin{enumerate}
\item for each $N\in\N$ the measure $\rho^{(N)}$ is probabilistic,
\item if for each $N\in\N$ we have that $A_1^{(N)}, A_2^{(N)}, A_3^{(N)}$ 
are independent random variables with distribution given by the measure 
$\rho^{(N)}$ then the probability of the event: $A_1^{(N)}\cap A_2^{(N)} 
\cap A_3^{(N)}\neq \emptyset$ tends to $0$ as $N$ tends to infinity,
\item let $k$ be a fixed natural number. if for each $N\in\N$ we have that
$A_1^{(N)},\dots,A_k^{(N)}$ are independent random variables with 
distribution given by the measure $\rho^{(N)}$ then the joint distribution 
of the random variables $|A_i\cap A_j|$, $1\leq i<j\leq k$ 
tends to a product distribution as $N$ tends to infinity. The limit 
distribution of a single variable $|A_i\cap A_j|$ is given by
$$\lim_{N\rightarrow\infty} P\big(|A_i^{(N)}\cap A_j^{(N)}|=k\big)= p_k.$$
\end{enumerate}
\end{remark}

\begin{remark}
The assumption \ref{zal:z4} follows from other assumptions, however the proof
of this is rather technical and we omit it.
\end{remark}

\begin{remark}
We would like to point out an interesting informal connection between 
stochastical properties of $G_i$ regarded as large matrices and their 
entries.

Let us consider $\Gamma_{\mu\nu}=\delta_{\mu\nu}$. 
For $q=1$ we have that $G_\mu$ is a family of independent Gaussian 
variables and for $q=0$ we have that $G_\mu$ is a family of free 
non commutative random variables. Freeness is an analogue of classical 
independence; we can expect therefore that for in a general case $-1\leq 
q\leq 1$ variables $G_\mu$ are independent in some generalised way.

However, it was proven by Speicher \cite{Sp3} that there are only three 
generalisations of the notion of independence of random variables 
to the non commutative setup which would satisfy certain natural properties.
These three generalisations are: classical independence, free independence 
(freeness) and boolean independence. Therefore except the cases 
$q\in\{0,1\}$ which correspond to the free and the classical situation 
respectively we cannot formally say that the non commutative random 
variables $G_\mu$ are independent in some sense.

We can of course weaken Speicher's axioms and treat this ``independence'' 
on a very informal level. It is worth pointing out that a family of 
``independent'' variables $G_i$ is asymptotically represented as
random matrices such that the entries of different matrices are classically
independent random variables.

Similarly, for the choice of the covariance function $\Gamma_{ts}=\min 
(t,s)$ for $t,s\geq 0$ we obtain a non commutative stochastic process $G_t$ 
which can be regarded as some kind of a Brownian motion \cite{BS0} and 
$G_t$ can be asymptotically represented as a matrix valued stochastic 
process. Entries of this matrix are classical Brownian motions.  
\end{remark}

\begin{remark} The assumptions of the theorem are fulfilled for the following
important examples of the functions $\sigma$: \end{remark}
\begin{proposition} \label{prop:zbieznosc}
For every real $c>0$ the sequence of functions defined by
$$(\sigma_A^{(N)})^2=\left( \frac{c}{\sqrt{N}} \right)^{|A|} \left(
1-\frac{c}{\sqrt{N}} \right)^{N-|A|}$$
fulfils the assumptions of Theorem \ref{th:theo} with the sequence 
$p_k=\frac{1}{k!} c^{2k} e^{-c^2}$ being the Poisson distribution with 
parameter $c^2$ and $$q=e^{-(1-\frac{1}{d^2}) c^2}.$$ 

In this case the covariance (\ref{eq:kowariancjas}) takes a beautiful form 
\begin{equation} \E[S^{\mu}_{\ii\jj} S^{\nu}_{\kk\lll}] =
\E[S^{\mu}_{\ii\jj} \overline{S^{\nu}_{\lll\kk}}]=
\label{eq:przykladkowariancji}
\end{equation}
$$= \Gamma_{\mu\nu} \prod_r \left(\frac{c}{\sqrt{N}} 
\frac{[i_r=l_r] [j_r=k_r]}{d} + 
\left(1-\frac{c}{\sqrt{N}}\right) [i_r=j_r] [k_r=l_r] \right).
$$
\end{proposition}
Proof of this proposition will be presented in Sect.\  
\ref{sec:dowodzbieznosci}. 
\begin{proposition}
For every real number $c>0$ the 
sequence of functions $\sigma^{(N)}$ defined for $N$ sufficiently large by 
$$(\sigma_A^{(N)})^2=\left\{ \begin{array}{cl} \frac{1}{{N \choose \lfloor c 
\sqrt{N} \rfloor }} &\quad \mbox{if } |A|=\lfloor c\sqrt{N}\rfloor \\ 0 & 
\quad \mbox{otherwise} \end{array} \right. , $$
where $\lfloor x \rfloor$ denotes the integer part of a real number $x$,
fulfils the assumptions of Theorem \ref{th:theo} with
$p_k=\frac{1}{k!} c^{2k} e^{-c^2}$ and $$q=e^{-(1-\frac{1}{d^2}) c^2}.$$
\end{proposition}
Since proof of this proposition is similar to the proof of Proposition 
\ref{prop:zbieznosc} we skip it.

\begin{lemma} \label{lem:glownylemat}
For any pair partition $\pi=\bigl\{ \{c_1,d_1\},\dots,\{c_m,d_m\} \bigr\}$ of
the set $\{1,\dots,2m\}$ and any sets 
$A_1,\dots,A_{2m}\subseteq\{1,\dots,N\}$ we have $$0\leq \frac{1}{d^N} 
\sum_{\ii^1,\dots,\ii^{2m}}\  \prod_{1\leq v\leq m} 
\T^{A_{c_v}}_{\ii^{c_v}\ii^{c_v+1}, \ii^{d_v}\ii^{d_v+1} } \leq 1.$$

If furthermore $A_{c_i} \cap A_{c_j} \cap A_{c_k}=\emptyset$ for all
$1\leq i<j<k\leq m$ then
$$\frac{1}{d^N} \sum_{\ii^1,\dots,\ii^{2m}}\  \prod_{1\leq v\leq m}
\T^{A_{c_v}}_{\ii^{c_v}\ii^{c_v+1},\ii^{d_v}\ii^{d_v+1}}=\prod_{{1\leq
i<j \leq m}\atop {{\rm lines\ }
\{c_i,d_i\}, \{c_j,d_j\}\ {\rm cross}}} \frac{1}{d^{2|A_{c_i}\cap A_{c_j}|}
}.$$
\end{lemma}
The proof of this lemma will be presented in Sect.\ 
\ref{sec:dowodzbieznosci}.

The following lemma states a well known property of the Gaussian distribution. 
\begin{lemma} \label{lem:gaussy}
If the joint distribution of random
variables $(X_k)$ is Gaussian and $\E[X_k]=0$ then
$$\E (X_1\cdots X_{2m-1})=0,$$
$$\E(X_1\cdots X_{2m})=\sum_\pi \E(X_{c_1} X_{d_1}) \cdots \E(X_{c_m}
X_{d_m}),$$
where the sum is taken over all pair partitions $\pi=\bigl\{
\{c_1,d_1\},\dots,\{c_m,d_m\} \bigr\}$ of
the set $\{1,\dots,2m\}$.
\end{lemma}

With this preparation we are able to start the proof of the main theorem.
\begin{proof}[Proof of Theorem \ref{th:theo}]
In the following the sums over $\pi$ are taken over all pair partitions
$\pi=\bigl\{ \{c_1,d_1\},\dots,\{c_m,d_m\} \bigr\}$ of the set
$\{1,\dots,2m\}$ and sums over $A_1,\dots,A_{n}$ are taken over all subsets
of the set $\{1,\dots,N\}$. Products over $v$ are taken over $1\leq v\leq m$. 

From Lemma \ref{lem:gaussy} follows that for any $m\in\N$
and indexes $\mu_1,\dots,\mu_{2m}\in\M$ we 
have: 
$$\E\ \tr(S^{(N),\mu_1} \cdots S^{(N),\mu_{2m-1}})= 0$$
and furthermore
$$ U^{(N)}:=\E\ \tr(S^{(N),\mu_1} \cdots S^{(N),\mu_{2m}})= $$
$$=\sum_\pi
\sum_{A_1,\dots,A_{2m} } \left[
\frac{\sigma^{(N)}_{A_1} \cdots \sigma^{(N)}_{A_{2m}} }{d^N}
\sum_{\ii^1,\dots,\ii^{2m}}  \prod_{v} 
\delta_{A_{c_v} A_{d_v}} \Gamma_{\mu_{c_v} \mu_{d_v}}
\T^{A_{c_v}}_{\ii^{c_v}\ii^{c_v+1},\ii^{d_v}\ii^{d_v+1}}
\right].$$

We define
$$V^{(N)}:=\sum_\pi \sum_{A_1,\dots,A_{2m} }
\Bigg[ \left(\prod_{v} \equalA{A_{c_v}}{A_{d_v} } \Gamma_{\mu_{c_v}
\mu_{d_v}} \right) \sigma^{(N)}_{A_1} \cdots \sigma^{(N)}_{A_{2m}} \times
$$ $$\times \prod_{{1\leq i<j \leq m}\atop {{\rm lines\ }
\{c_i,d_i\}, \{c_j,d_j\}\ {\rm cross}}} \frac{1}{d^{2|A_{c_i}\cap A_{c_j}|}
} \Bigg]$$

Lemma \ref{lem:glownylemat} shows that the corresponding summands in
the definitions of $U$ and $V$ are equal unless there are some
indexes $1\leq p<q<r\leq m$ such that $A_{c_p}\cap A_{c_q} \cap
A_{c_r}\neq\emptyset$. There are $m \choose 3$ choices of these indexes and
again from Lemma \ref{lem:glownylemat} and the assumption \ref{zal:z1} we 
have 
$$|U^{(N)}-V^{(N)}|\leq C {m \choose
3} \sum_{{A_1,A_2,A_3\subseteq\{1,\dots,N\}}\atop{A_1\cap A_2\cap
A_3\neq\emptyset}} \ \sum_{A_4,\dots,A_m\subseteq\{1,\dots,N\} }
(\sigma^{(N)}_{A_1})^2 \cdots (\sigma^{(N)}_{A_m})^2=$$
$$=C {m \choose 3} \sum_{{A_1,A_2,A_3\subseteq\{1,\dots,N\}}\atop{A_1\cap 
A_2\cap A_3\neq\emptyset}} (\sigma^{(N)}_{A_1})^2 (\sigma^{(N)}_{A_2})^2
(\sigma^{(N)}_{A_3})^2, $$ 
where $C=\max \left| \Gamma_{\mu_p,\mu_q} \right|$ 
and therefore from the assumption \ref{zal:z2} we 
have $$\lim_{N\rightarrow\infty} |U^{(N)}-V^{(N)}|=0.$$

We have
$$V^{(N)}=\sum_\pi  \left( \prod_v \Gamma_{\mu_{c_v} \mu_{d_v}} \right)
 \sum_{A_1,\dots,A_m} (\sigma^{(N)}_{A_1})^2 \cdots
(\sigma^{(N)}_{A_m})^2 \times $$ $$\times \prod_{{1\leq i<j \leq m}\atop 
{{\rm lines\ } \{c_i,d_i\}, \{c_j,d_j\}\ {\rm cross}}} 
\frac{1}{d^{2|A_{c_i}\cap A_{c_j}|} }=$$
$$=\sum_\pi  \left( \prod_v \Gamma_{\mu_{c_v} \mu_{d_v}} \right)
\int \prod_{{1\leq i<j \leq m}\atop 
{{\rm lines\ } \{c_i,d_i\}, \{c_j,d_j\}\ {\rm cross}}} 
\frac{1}{d^{2n_{ij}} } d\lambda^{(N)}(n_{ij}),$$
%
where measures $(\lambda^{(N)})$ are defined on the set of all sequences 
$(n_{ij})_{1\leq i<j\leq m}$, $n_{ij}\in\{0,1,2,\dots\}$ by condition 
$$\lambda^{(N)} 
\big(\{(n_{ij})\}\big)=\sum_{{A_1,\dots,A_k\subseteq\{1,\dots,N\}} \atop 
{|A_i\cap A_j|=n_{ij} {\rm\ for\ all\ } 1\leq i<j\leq m}} 
(\sigma^{(N)}_{A_1})^2 \cdots (\sigma^{(N)}_{A_N})^2. $$ 
From the 
assumption \ref{zal:z3} follows that this sequence converges 
pointwise to the product 
measure defined on the atoms by
$$\lambda \big(\{(n_{ij}) \}\big)=\prod_{1\leq i<j\leq m} p_{ij}.$$

Since measures $\lambda^{(N)}$ and the measure $\lambda$ are probabilistic, 
this convergence is uniform and the statement of the theorem follows.
\end{proof}

\section{The Almost Surely Convergence}
\label{sec:almostsurely}
As a simple corollary of Lemma \ref{lem:gaussy} we have
\begin{lemma} \label{lem:smiesznesumowanie}
If the joint distribution of random variables $X_1,\dots,X_{2m}$ is Gaussian
and $\E[X_k]=0$ for all $1\leq k\leq 2m$ then
$$\E[X_1 \cdots X_{2m}]-\E[X_1 \cdots X_m]\ \E[X_{m+1} \cdots X_{2m}]=$$
$$=\sum_{\pi'} \prod_{1\leq v\leq m} \E[X_{c_v} X_{d_v} ],$$
where the sum is taken over pair partitions
$\pi'=\big\{ \{c_1,d_1\},\dots,\{c_m,d_m\} \big\}$ of the set
$\{1,\dots,2m\}$ which have additional property that there exist
$x\in\{1,\dots,m\}$ and $y\in\{m+1,\dots,2m\}$ such that $\{x,y\}\in\pi'$.
\end{lemma}
We define a permutation $\sigma:\{1,\dots,2m\}\rightarrow\{1,\dots,2m\}$ by
$\sigma(k)=k+1$ for $k\not\in\{m,2m\}$, $\sigma(m)=1$, $\sigma(2m)=m+1$.
\begin{lemma}
Let $\pi'=\big\{ \{c_1,d_1\},\dots,\{c_m,d_m\} \big\}$ be a pair partition 
as in Lemma \ref{lem:smiesznesumowanie}, i.e. for some $1\leq k\leq m$ 
we have $c_k\in\{1,\dots,m\}$, $d_k\in\{m+1,\dots,2m\}$.

Then 
$$\left| \frac{1}{d^{2N}} \sum_{\ii^1,\dots,\ii^{2m}} \prod_{1\leq v\leq m} 
\T^{A_{c_v}}_{\ii^{c_v} \ii^{\sigma(c_v)},\ii^{d_v} \ii^{\sigma(d_v)}} 
\right|\leq \frac{1}{d^{2 |A_k\setminus \bigcup_{j\neq k} A_j |} }
$$
\end{lemma}
This lemma follows directly from Lemmas \ref{lem:produktowo} and 
\ref{lem:pojedynczyiloczyn2} from Sect.\ \ref{sec:dowodzbieznosci}.

\begin{proposition} \label{prop:wariancja}
$${\rm Var}[\tr\ S^{(N),\mu_1} \cdots S^{(N),\mu_m}]\leq  $$
$$\leq C \sum_{A_1,\dots,A_m\subseteq\{1,\dots,N\} } (\sigma^{(N)}_{A_1})^2 
\cdots (\sigma^{(N)}_{A_m})^2 \frac{1}{d^{2|A_1\setminus (A_2\cup 
\cdots\cup A_m)|}}, $$
where $C=(2m)!! \max |\Gamma_{\mu_p,\mu_q}|$. 
\end{proposition} 
\begin{proof}
We define $\mu_{m+k}=\mu_k$. In the following sums over $\pi´$ are taken over all
pair partitions $\pi´=\bigl\{ \{c_1,d_1\},\dots,\{c_m,d_m\} \bigr\}$ of the set
$\{1,\dots,2m\}$ with additional property that there exist
$x\in\{1,\dots,m\}$ and $y\in\{m+1,\dots,2m\}$ such that $\{x,y\}\in\pi'$.

From Lemma \ref{lem:smiesznesumowanie} we 
have that $${\rm Var}[\tr\ S^{(N),\mu_1} 
\cdots S^{(N),\mu_m}]=$$ $$=\E [(\tr\ S^{(N),\mu_1} \cdots S^{(N),\mu_m})^2]-
[\E(\tr\ S^{(N),\mu_1} \cdots S^{(N),\mu_m})]^2=$$
$$=\frac{1}{d^{2N}} \sum_{\ii^1,\dots,\ii^{2m} }
\bigg( \E[S^{\mu_1}_{\ii^1 \ii^2} \cdots S^{\mu_{m-1}}_{\ii^{m-1} \ii^m} 
S^{\mu_m}_{\ii^m \ii^{1}}\ S^{\mu_1}_{\ii^{m+1} \ii^{m+2}} \cdots 
S^{\mu_{m-1}}_{\ii^{2m-1} \ii^{2m}} S^{\mu_m}_{\ii^{2m} \ii^{m+1}} ]-$$
$$-\E[S^{\mu_1}_{\ii^1 \ii^2} \cdots S^{\mu_{m-1}}_{\ii^{m-1} \ii^m} 
S^{\mu_m}_{\ii^m \ii^{1}}] \ \E[S^{\mu_1}_{\ii^{m+1} \ii^{m+2}}
\cdots S^{\mu_{m-1}}_{\ii^{2m-1} \ii^{2m}} S^{\mu_m}_{\ii^{2m} \ii^{m+1}} ] 
\bigg)=$$ 
$$=\frac{1}{d^{2N}} \sum_{\pi'} \sum_{\ii^1,\dots,\ii^{2m}} 
\prod_{1\leq v\leq m} \E[S^{\mu_{c_v}}_{\ii^{c_v} \ii^{\sigma(c_v)}} 
S^{\mu_{d_v}}_{\ii^{d_v} \ii^{\sigma(d_v)}} ]=$$
$$=\frac{1}{d^{2N}} \sum_{\pi'} 
\sum_{A_1,\dots,A_{2m}\subseteq\{1,\dots,N\} } \sum_{\ii^1,\dots,\ii^{2m}}
\sigma^{(N)}_{A_1} \cdots \sigma^{(N)}_{A_{2m}}\times$$
$$\times \prod_{1\leq v\leq m} 
\equalA{A_{c_v}}{A_{d_v}} \Gamma_{\mu_{c_v},\mu_{d_v}} 
\T^{A_{c_v}}_{\ii^{c_v} \ii^{\sigma(c_v)},\ii^{d_v} \ii^{\sigma(d_v)}}\leq$$
$$\leq C \sum_{A_1,\dots,A_m\subseteq\{1,\dots,N\} }
(\sigma^{(N)}_{A_1})^2 \cdots (\sigma^{(N)}_{A_m})^2 
\frac{1}{d^{2|A_1\setminus (A_2\cup \cdots\cup A_m)|}}, $$
where in the last inequality we used Lemmas \ref{lem:produktowo} 
and \ref{lem:pojedynczyiloczyn2}.
\end{proof}

If $(M_{ij})_{1\leq i,j\leq K}$ is a hermitian matrix with eigenvalues 
$\lambda_1,\dots,\lambda_K$ we define a probabilistic measure $\nu_M$ on 
the real line $\R$ by the following 
$$\nu_M=\frac{1}{K} \sum_{1\leq n\leq K} \delta_{\lambda_n}.$$

\begin{theorem}
If the assumptions of Theorem \ref{th:theo} are fulfilled then there 
exists an increasing sequence of natural numbers $(N_i)$ such that 
 the sequence of measures $\nu_{S^{(N_i)}}$ almost 
surely converges weakly to the measure $\nu_q$. 
\end{theorem}

Since for $0\leq q<1$ the support of the limit measure $\nu_q$ is compactly 
supported, the theorem follows from the following stronger statement. 

\begin{theorem}
If the assumptions of Theorem \ref{th:theo} are fulfilled then there 
exists a sequence $(N_i)$ such that the limit
$$ \lim_{m\rightarrow\infty} \tr\ S^{(N_i),\mu_1} \cdots S^{(N_i),\mu_m} =  
\tau(G_{\mu_1} \cdots G_{\mu_m}) $$
holds almost surely.
\end{theorem}
\begin{proof}
Our goal is to construct a sequence $(N_i)$ such that
$$\sum_{i}  \E\ \Big[\big( \tr\ S^{(N_i),\mu_1} \cdots S^{(N_i),\mu_m} - 
\tau(G_{\mu_1} \cdots G_{\mu_m}) \big)^2 \Big]<\infty $$
holds.

However,
$$ \E\ \Big[\big( \tr\ S^{(N),\mu_1} \cdots S^{(N),\mu_m} -
\tau(G_{\mu_1} \cdots G_{\mu_m})\big)^2 \Big]=$$ 
$$ \big( \E\ [\tr\ S^{(N),\mu_1} 
\cdots S^{(N),\mu_m}- \tau(G_{\mu_1} \cdots G_{\mu_m})]\big)^2 +
{\rm Var}\ 
[\tr\ S^{(N),\mu_1} \cdots S^{(N),\mu_m}]. $$
The first summand converges to $0$ by Theorem \ref{th:theo}.

From Proposition \ref{prop:wariancja} we have that
$${\rm Var}[\tr\ S^{(N),\mu_1} \cdots S^{(N),\mu_m}]\leq$$  
$$\leq C \sum_{A_1,\dots,A_m\subseteq\{1,\dots,N\} } (\sigma^{(N)}_{A_1})^2 
\cdots (\sigma^{(N)}_{A_m})^2 \frac{1}{d^{2|A_1\setminus (A_2\cup 
\cdots\cup A_m)|}}$$  
From assumptions \ref{zal:z1} and \ref{zal:z4} it follows that this 
expression converges to $0$. \end{proof}


\section{Proofs of Technical Lemmas}
\label{sec:dowodzbieznosci}
\label{sec:technicznedowody}
Lemma \ref{lem:glownylemat} follows directly from the following two lemmas.
\begin{lemma} \label{lem:produktowo}
For every $n,M\in\N$ if for all $1\leq v\leq M$ we have
$A_v\subseteq\{1,\dots,N\}$, $1\leq p_v,q_v,r_v,s_v\leq n$ then
\begin{equation} \label{eq:formulazlematu} 
\sum_{\ii^1,\dots,\ii^n}
\prod_{1\leq v\leq M} \T^{A_v}_{\ii^{p_v} \ii^{q_v},
\ii^{r_v} \ii^{s_v}} =
\prod_{1\leq r\leq N}\ \sum_{0\leq j^1,\dots,j^{n}\leq d-1}\
\prod_{1\leq v\leq M} T^{A_v,r}_{j^{p_v}   j^{q_v}  ,j^{r_v}   j^{s_v}  }.
\end{equation}
\end{lemma}
\begin{proof} This lemma is a direct consequence of 
Eq.\ (\ref{eq:tensort1}) and (\ref{eq:tensort2}). \end{proof}

For sets $A_1,\dots,A_{2m}\subseteq\{1,\dots,N\}$ and a pair partition
$\pi=\bigl\{ \{c_1,d_1\},\dots,\{c_m,d_m\} \bigr\}$ of the set
$\{1,\dots,2m\}$ we define \begin{equation}
\Theta^{A_1,\dots,A_{2m},\pi}_r=
\frac{1}{d} \sum_{0\leq j^1,\dots,j^{2m}\leq d-1} 
\prod_u T^{A_{c_u},r}_{j_{c_u} j_{c_u+1},j_{d_u} j_{d_u+1}}= 
\label{eq:slynnewyrazenie} \end{equation}
$$=\frac{1}{d}\sum_{0\leq j^1,\dots,j^{2m}\leq
d-1}\ \prod_{1\leq u\leq m} \left\{ \begin{array}{cc} 
\equal{j^{c_u}}{j^{c_u+1} } \equal{j^{d_u}}{j^{d_u+1} } &\quad {\rm if\ }  
r\not\in A_{c_u}    \\ 
\frac{1}{d} \equal{j^{c_u}}{j^{d_u1} } 
\equal{j^{c_u+1}}{j^{d_u} } &\quad {\rm if\ } r\in A_{c_u} \end{array} \right. . $$
It should be understood that $j^{2m+1}=j^1$.

\begin{lemma} \label{lem:pojedynczeczynniki1}
For any $A_1,\dots,A_{2m}\subseteq\{1,\dots,N\}$ and a pair partition
$\pi=\bigl\{ \{c_1,d_1\},\dots,\{c_m,d_m\} \bigr\}$ we have
\begin{equation} 0\leq \prod_{1\leq r\leq N}
\Theta^{A_1,\dots,A_{2m},\pi}_r\leq 1. \end{equation}
If furthermore $A_{c_i} \cap A_{c_j} \cap A_{c_k}=\emptyset$ for all
$1\leq i<j<k\leq n$ then
\begin{equation} \prod_{1\leq r\leq N}
\Theta^{A_1,\dots,A_{2m},\pi}_r=\prod_{{1\leq
i<j \leq m}\atop{{\rm lines\ }
\{c_i,d_i\}, \{c_j,d_j\}\ {\rm cross}}} \frac{1}{d^{2|A_{c_i}\cap A_{c_j}|}
}.\end{equation}

\end{lemma}

\begin{proof} 
Let us consider an expression of the following type
\begin{equation} \sum_{0\leq j^1,\dots,j^{2m}\leq d} [j^{e_1}=j^{f_1}] \cdots [j^{e_{2m}}=j^{f_{2m}}]. 
\label{eq:sumaiversona}
\end{equation}
We can represent this expression by a graph (see for example \cite{Do}) 
with $2m$ vertices which are labelled by variables
$j^1,\dots, j^{2m}$ and with vertices $j^{e_i}$ and $j^{f_i}$ connected by an edge for all $1\leq i\leq 2m$.
We see that the nonzero summands in (\ref{eq:sumaiversona}) come from indexes $(j^1,\dots,j^{2m})$ such that
to all vertices in the same connected component of the graph is assigned the same value.
Therefore the expression (\ref{eq:sumaiversona}) is equal to $d^M$, where $M$ denotes the number of
connected components of the graph. 

Let us fix an index $r$. We shall apply the above observation to 
evaluate $\Theta^{A_1,\dots,A_{2m},\pi}_r$.
Let $n$ be the number of indexes $v$ such that $r\in A_{c_v}$ and let 
$v_1,\dots,v_n$ be all such indexes. We consider the following cases.
\begin{enumerate}
\item If $n=0$ then we obtain a graph of type presented in Fig.\ \ref{fig:fig1}. This graph has one
connected component and
therefore $\Theta_r^{A_1,\dots,A_{2m},\pi}=1$.
\item If $n=1$ then we obtain a graph of type presented in Fig.\ \ref{fig:fig2}. This graph has two
connected components and therefore
$\Theta^{A_1,\dots,A_{2m},\pi}_r=1$.
\item Suppose that $n=2$. If the lines $\{c_{v_1},d_{v_1}\}$ and
$\{c_{v_2},d_{v_2}\}$ do not cross (Fig.\ \ref{fig:fig3} and \ref{fig:fig4}) then corresponding
graphs have three components and and therefore
$\Theta^{A_1,\dots,A_{2m},\pi}_r=1$.

But if the lines $\{c_{v_1},d_{v_1}\}$
and $\{c_{v_2},d_{v_2}\}$ cross (Fig.\ \ref{fig:fig5}) then the graph has 
only one connected  component and
therefore $\Theta^{A_1,\dots,A_n,\pi}_r=\frac{1}{d^2}$.

\item In the general case $n\geq 3$ the graph has at most $n+1$ components and
there is a factor $\frac{1}{d^{n+1}}$, therefore
$0<\Theta^{A_1,\dots,A_{2m},\pi}_r\leq 1$.
\end{enumerate}

\begin{figure}  
\centering
\input{new1.pstex_t}
\caption{}
\label{fig:fig1}
\end{figure}
\begin{figure}  
\centering
\input{new2.pstex_t}
\caption{}
\label{fig:fig2}
\end{figure}
\begin{figure}  
\centering
\input{new3.pstex_t}
\caption{}
\label{fig:fig3}
\end{figure}
\begin{figure}  
\centering
\input{new4.pstex_t}
\caption{}
\label{fig:fig4}
\end{figure}
\begin{figure}  
\centering
\input{new5.pstex_t}
\caption{}
\label{fig:fig5}
\end{figure}

Therefore if for a given partition $\pi$ sets $A_{c_1},\dots,A_{c_k}$ have
a property that the common part of each three of them is empty then
\begin{equation} \prod_{1\leq r\leq N}
\Theta^{A_1,\dots,A_n,\pi}_r=\prod_{{1\leq
i<j \leq k}\atop{{\rm lines\ }
\{c_i,d_i\}, \{c_j,d_j\}\ {\rm cross}}} \frac{1}{d^{2 |A_{c_i}\cap A_{c_j}|}
}.\end{equation}

In the general situation the expression $\prod_{1\leq r\leq N}
\Theta^{A_1,\dots,A_{2m},\pi}_r$ is a real number
from the interval $[0,1]$.
\end{proof}


We recall that the 
permutation $\sigma:\{1,\dots,2m\}\rightarrow\{1,\dots,2m\}$ was defined by
$\sigma(k)=k+1$ for $k\not\in\{m,2m\}$, $\sigma(m)=1$, $\sigma(2m)=m+1$.

For sets $A_1,\dots,A_{2m}\subseteq\{1,\dots,N\}$ and a pair partition
$\pi=\bigl\{ \{c_1,d_1\},\dots,\{c_m,d_m\} \bigr\}$ of the set
$\{1,\dots,2m\}$ we define 
\begin{equation} 
\Upsilon^{A_1,\dots,A_{2m},\pi}_r=\frac{1}{d}\sum_{0\leq j^1,\dots,j^{2m}\leq
d-1}\ \prod_{1\leq v\leq m} T^{A_c,r}_{c_v \sigma(c_v),d_v \sigma(d_v)}=
 \label{ref:slynnewyrazenie7} \end{equation}
$$=\frac{1}{d}\sum_{0\leq j^1,\dots,j^{2m}\leq
d-1}\ \prod_{1\leq v\leq m} \left\{ \begin{array}{cc} 
\equal{j^{c_v}}{j^{\sigma(c_v)} } \equal{j^{d_v}}{j^{\sigma(d_v)} } & {\rm 
if\ }  r\not\in A_{c_v} 
\\ \frac{1}{d} \equal{j^{c_v}}{j^{\sigma(d_v)} } 
\equal{j^{\sigma(c_v)}}{j^{d_v} } & {\rm if\ } r\in A_{c_v} \end{array} 
\right. .$$ 
\begin{lemma} \label{lem:pojedynczyiloczyn2}
For any $A_1,\dots,A_{2m}\subseteq\{1,\dots,N\}$ and a pair partition
$\pi=\bigl\{ \{c_1,d_1\},\dots,\{c_m,d_m\} \bigr\}$ we have
\begin{equation} 0\leq \Upsilon^{A_1,\dots,A_{2m},\pi}_r\leq 1. 
\end{equation} If furthermore $r\in A_k\setminus\bigcup_{j\neq k} A_j$ then
\begin{equation} 
\Upsilon^{A_1,\dots,A_{2m},\pi}_r=\frac{1}{d^2}
\end{equation} \end{lemma} 
Proof of this lemma is very similar to the proof of Lemma 
\ref{lem:pojedynczeczynniki1} and we will skip it. 

\begin{proof}[Proof of Prop.\ \ref{prop:zbieznosc}]
The measure $\rho^{(N)}$ on the set of all subsets of 
$\{1,\dots,N\}$ such that for any $A\subseteq\{1,\dots,N\}$ we have
$$\rho^{(N)}(\{A\})=\left( \frac{c}{\sqrt{N}} \right)^{|A|} \left(
1-\frac{c}{\sqrt{N}} \right)^{N-|A|}$$ is simply a Bernoulli distribution and 
the assumption \ref{zal:z1} of Theorem \ref{th:theo} is obviously 
fulfilled.

Therefore if $A_1,A_2,A_3$ are independent random variables with 
distribution given by $\rho^{(N)}$ then
$$P\{\omega:A_1(\omega)\cap A_2(\omega) \cap A_3(\omega)\neq\emptyset\}\leq$$
$$\leq\sum_{1\leq k\leq N} P\{\omega:k\in A_1(\omega)\cap A_2(\omega) \cap 
A_3(\omega)\}=N \left(\frac{c}{\sqrt{N} } \right)^3 $$
and the assumption \ref{zal:z2} of Theorem \ref{th:theo} is fulfilled.

Let $A_1,\dots,A_k$ be independent random variables with distribution given 
by $\rho^{(N)}$ and let us consider the probability of following event: 
$|A_p\cap 
A_q|=n_{pq}$. The discussion from the previous paragraph shows that this 
probability
is equal (up to an error of order $O(N^{-\frac{1}{2}})$)  to the probability 
of the 
event: $|A_p\cap A_q|=n_{pq}$ and furthermore $A_p\cap A_q \cap A_r=\emptyset$
for all $p<q<r$. We shall now evaluate the latter probability.

There are $\frac{N!}{\left( \prod n_{pq}! \right) \left(N-\sum 
n_{pq} \right)!}$ choices of disjoint sets $(B_{pq})_{1\leq p<q\leq k}$ such 
that $|B_{pq}|=n_{pq}$. 

For $i\in B_{pq}$ the probability of the event: $i\in A_{p}\cap A_{q}$ and 
$i\not\in A_r$ for $r\not\in\{p,q\}$ is equal to
$$\left(\frac{c}{\sqrt{N}}\right)^2 
\left(1-\frac{c}{\sqrt{N}}\right)^{k-2}.$$

For $i\not\in \bigcup_{pq} B_{pq}$ the probability of the event: $i\not\in 
A_{p}\cap A_{q}$ for all $1\leq p<q\leq k$ is equal to
$$\left(1-\frac{c}{\sqrt{N}} \right)^k+k \frac{c}{\sqrt{N}} 
\left(1-\frac{c}{\sqrt{N}}\right)^{k-1}.$$

Therefore the probability of the considered event is equal to
$$\frac{N!}{\left( \prod n_{pq}! \right) 
\left(N-\sum n_{pq} \right)!}
\left(\frac{c}{\sqrt{N} }\right)^{2\sum n_{pq}} 
\left(1-\frac{c}{\sqrt{N}}\right)^{(k-2) \sum n_{pq}} \times $$ $$\times
\left[\left(1-\frac{c}{\sqrt{N}} \right)^k+k \frac{c}{\sqrt{N}} 
\left(1-\frac{c}{\sqrt{N}}\right)^{k-1} \right]^{N-\sum n_{pq}} $$
where all sums and products are taken over $1\leq p<q\leq k$.

After short calculations one can show that the limit of this expression as
$N$ tends to infinity is equal to
$$\prod_{1\leq p<q\leq k} \frac{1}{n_{pq}!} c^{2n_{pq}} e^{-c^2} $$
and therefore the assumption \ref{zal:z3} is fulfilled.

The assumption \ref{zal:z4} follows from the observation that the 
distribution of the random variable 
$|A_1\setminus(A_2\cup\cdots A_n)|$ is binomial and simple computations.
\end{proof}

\section{Acknowledgements}
I would like to thank Roland Speicher and Marek Bo\.zejko for very 
inspiring discussions and Franz Lehenr for numerous remarks concerning the 
manuscript.


The work was partially supported by the Scientific Committee in Warsaw under
grant number P03A05415. The paper was written while the author was on leave
in the University of Heidelberg on a scholarship funded by German Academic
Exchange Service (DAAD).


\begin{thebibliography}{999}
\bibitem [Boz]{B} Bo\.zejko, M:
Completely positive maps on Coxeter groups and the ultracontractivity of the 
$q$-Ornstein-Uhlenbeck semigroup.
Alicki, Robert (ed.) et al., Quantum probability. Workshop, Gdansk, Poland, July 
1--6, 1997. Warsaw: Polish
Academy of Sciences, Institute of Mathematics, Banach Cent. Publ. \textbf{43}, 87-93 
(1998) 
\bibitem [BKS]{BKS} Bozejko, M., Kuemmerer, B.\ and Speicher, R.: $q$-Gaussian 
processes: Non-commutative and classical aspects. 
Commun.\ Math.\ Phys.\ \textbf{185}, 129--154 (1997)
\bibitem [BS1]{BS0} Bo\.zejko, M.\ and Speicher, R.:
An example of a generalized Brownian motion.
Commun. Math. Phys. \textbf{137}, 519--531 (1991).
\bibitem [BS2]{BS} Bo\.zejko, M.\ and Speicher, R.: Completely positive maps on Coxeter
groups, deformed commutation relations and operator spaces, Math. Ann. \textbf{300}, 97--120
(1994)  
\bibitem [Br]{Br} Brody, T.A., et al: Random--matrix physics: spectrum and strength fluctuations. 
Rev.\ Mod.\ Phys.\ \textbf{53}, 385--479 (1981)
\bibitem [Do]{Do} Douglas, M.R.:Large $N$ quantum field theory and matrix 
models. Voiculescu, Dan-Virgil (ed.), Free probability theory. Papers from 
a workshop on random matrices and operator algebra free products, Toronto, 
Canada, Mars 1995. Providence: American Mathematical Society. Fields Inst. 
Commun. \textbf{12}, 21-40 (1997)
\bibitem [FB]{FB} Frisch, U.\ and Bourret, R.: Parastochastics. J.\ Math.\ Phys. \textbf{11},
364--390 (1970)
\bibitem [Gr]{Greenberg} Greenberg, O.W.: Particles with small violations of Fermi or
Bose statistics. Phys.\ Rev.\ D \textbf{43}, 4111--4120 (1991)
\bibitem [GKP]{GKP} Graham, R.L., Knuth, D.E.\ and Patashnik, O.:
Concrete mathematics: a foundation for computer science. 2nd ed.
Amsterdam: Addison-Wesley Publishing Group, 1994
\bibitem [GMGW]{GMGW} Guhr, T., M\"uller-Groeling, A.\ and Weidenmüller, H.A.:
Random Matrix Theories in Quantum Physics: Common Concepts.  Phys. Rep. \textbf{299}, 
190--425 (1998) 
\bibitem [HP]{HP} Hiai, F.\ and Petz, D.: The semicircle law, free random variables
and entropy, to be published by AMS 
\bibitem [HT]{HT} Haagerup, U.\ and Thorbjoernsen, S.:
Random matrices and $K$-theory for exact $C^*$-algebras. Doc.\ Math., J.\ DMV \textbf{4}, 
341--450 (1999)
\bibitem [M]{M} Mehta, M.L.:
Random matrices. Rev.\ and enlarged 2.\ ed. Boston: Academic Press, 1991
\bibitem [MN]{MN} Mingo, J.\ and Nica, A.: Random unitaries in non--commutative
tori, and an asymptotic model for $q$--circular systems. Preprint.
\bibitem [Sh]{Sh} Shlyakhtenko, D.:
Random Gaussian band matrices and freeness with amalgamation.
Int.\ Math.\ Res.\ Not. 1996, No.\ 20, 1013-1025 (1996).
\bibitem [Sn1]{Sn1} \'Sniady, P.: On $q$--deformed quantum stochastic calculus. 
University of Wroclaw MSc thesis, Preprint 1999
\bibitem [Sp1]{Sp} Speicher, R.: A non--commutative central limit theorem. Math.\ Z.\
\textbf{209}, 55--66 (1992)
\bibitem [Sp2]{Speicher} Speicher, R.: Generalized Statistics of Macroscopic Fields.
Lett.\ Math.\ Phys.\ \textbf{27}, 97--104 (1993)
\bibitem [Sp3]{Sp3} Speicher, R.: On universal products.
Voiculescu, Dan-Virgil (ed.), Free probability theory. Papers from 
a workshop on random matrices and operator algebra free products, Toronto, 
Canada, Mars 1995. Providence: American Mathematical Society. Fields Inst. 
Commun. \textbf{12}, (1997)   
\bibitem [Sz]{Sz} Szego, G.: Ein Betrag zur Theorie der Thetafunktionen. Sitz.
Preuss. Akad. Wiss. Phys. Math. L1 19, 242--252 (1926)
\bibitem [V1]{V2} Voiculescu, D.:
Noncommutative random variables and spectral problems in free product 
$C\sp*$-algebras. Rocky Mt.\ J.\ Math.\ \textbf{20}, 263-283 (1990)
\bibitem [V2]{V3} Voiculescu, D.:
Limit laws for random matrices and free products. 
Invent.\ Math.\ \textbf{104}, 
201--220 (1991)
\bibitem [V3]{V4} Voiculescu, D.:
Free probability theory: Random matrices and von Neumann algebras. 
Chatterji, S. D. (ed.), Proceedings of the international congress of 
mathematicians, ICM '94, August 3-11,
1994, Zuerich, Switzerland. Vol. I. Basel: Birkhaeuser. 227-241 (1995). 
\bibitem [VDN]{VDN} Voiculescu, D.V., Dykema, K.J.\ and Nica, A.:
Free random variables. A noncommutative probability approach to free products 
with applications to
random matrices, operator algebras and harmonic analysis on free groups. CRM 
Monograph Series. 1. 
Providence: American Mathematical Society.,  1992

\end{thebibliography}
\end{document}